\documentclass[12pt]{article}
\usepackage{graphicx}
\usepackage{amsmath}
\usepackage{amsfonts}
\usepackage{amssymb}
\newtheorem{theorem}{Theorem}

\newtheorem{lemma}[theorem]{Lemma}

\newtheorem{proposition}[theorem]{Proposition}

\newenvironment{proof}[1][Proof]{\textbf{#1.} }{\ \rule{0.5em}{0.5em}}
\def\text{\hbox} \def\newpage{\vfill\break}

\def\a{\alpha}
\def\b{\beta}
\def\bb{\gamma}
\def\aa{\eta}

\def\la{\alpha}

\def\lbb{\gamma}

\def\e{\varepsilon}

\def\p{\pi}

\def\r{\rho}

\def\t{\textit{t}}

\def\ps{\varphi}

\def\f{\psi}

\def\o{\omega}

\def\O{\Omega}
\def\S{{\bf S}}

\def\T{{\cal T}}

\def\A{{\cal A}}

\def\K{{\cal K}}
\def\A'{{\cal A}'}
\def\Z{{\mathbb Z}}

\date{}

\begin{document}

\title{(1,1)-knots via the mapping class group of the twice punctured torus}

\author{Alessia Cattabriga \and Michele Mulazzani}

\maketitle


\begin{abstract}
We develop an algebraic representation for $(1,1)$-knots using the
mapping class group of the twice punctured torus $MCG_2(T)$. We
prove that every $(1,1)$-knot in a lens space $L(p,q)$ can be
represented by the composition of an element of a certain rank two
free subgroup of $MCG_2(T)$ with a standard element only depending
on the ambient space. As notable examples, we obtain a
representation of this type for all torus knots and for all
two-bridge knots. Moreover, we give explicit cyclic presentations
for the fundamental groups of the cyclic branched coverings of
torus knots of type $(k,ck+2)$.\\
\\ {{\it Mathematics Subject
Classification 2000:} Primary 57M05, 20F38; Secondary 57M12, 57M25.}

\end{abstract}


\section{Introduction and preliminaries} \label{intro}

The topological properties of $(1,1)$-knots, also called genus one
1-bridge knots, have recently been investigated in several papers
(see \cite{Be, CM, CK, Fu, Ga, Ga2, GM, Ha, Ha2, Ha3, MS, MSY,
MSY2, Mu, W1, Wu, W2}). These knots are very important in the
light of some results and conjectures involving Dehn surgery on
knots (see in particular \cite{Ga} and \cite{Wu}). Moreover, the
strict connection between cyclic branched coverings of
$(1,1)$-knots and cyclic presentations of groups have been pointed
out in \cite{CM}, \cite{GM} and \cite{Mu}.

Roughly speaking, a $(1,1)$-knot is a knot which can be obtained
by gluing along the boundary two solid tori with a trivial arc
properly embedded. A more formal definition follows. A set of
mutually disjoint arcs $\{t_1,\ldots ,t_b\}$ properly embedded in
a handlebody $H$ is \textit{trivial} if there exist $b$ mutually
disjoint discs $D_1,\ldots ,D_b\subset H$ such that $t_i\cap D_i=
t_i\cap\partial D_i=t_i$, $t_i\cap D_j=\emptyset$ and $\partial
D_i-t_i\subset\partial H$ for all $i,j=1,\ldots ,b$ and $i\neq j$.
Let $M=H\cup_{\varphi} H'$ be a genus $g$ Heegaard splitting of a
closed orientable 3-manifold $M$ and let $F=\partial H=\partial
H'$; a link $L\subset M$ is said to be in {\it $b$-bridge
position\/} with respect to $F$ if: (i) $L$ intersects $F$
transversally and (ii) $L\cap H$ and $L\cap H'$ are both the union
of $b$ mutually disjoint properly embedded trivial arcs. The
splitting is called a {\it $(g,b)$-decomposition of $L$\/}. A link
$L$ is called a {\it $(g,b)$-link\/} if it admits a
$(g,b)$-decomposition. Note that a $(0,b)$-link is a link in
$\S^3$ which admits a $b$-bridge presentation in the usual sense.
So the notion of $(g,b)$-decomposition of links in 3-manifolds
generalizes the classical bridge (or plat) decomposition of links
in $\S^3$ (see \cite{Do}). Obviously, a $(g,1)$-link is a knot,
for every $g\geq 0$.

\begin{figure}[ht]
\begin{center}
\includegraphics*[totalheight=7cm]{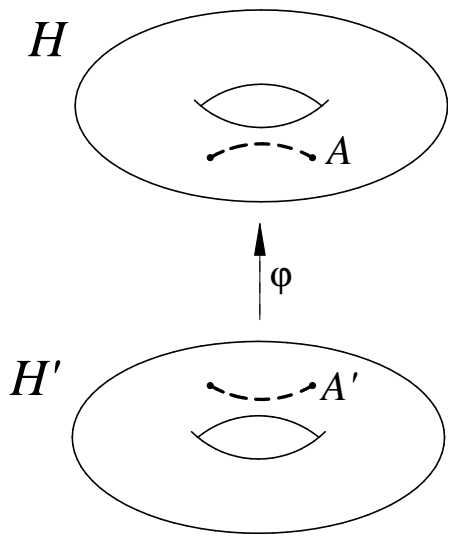}
\end{center}
\caption{A $(1,1)$-decomposition.} \label{Fig. 1}
\end{figure}

Therefore, a {\it $(1,1)$-knot\/} $K$ is a knot in a lens space
$L(p,q)$ (possibly in $\S^3$) which admits a $(1,1)$-decomposition
$$(L(p,q),K)=(H,A)\cup_{\ps}(H',A'),$$
where $\ps:(\partial H',\partial A')\to(\partial H,\partial A)$ is
an (attaching) homeomorphism which reverses the standard orientation
on the tori (see Figure \ref{Fig. 1}). It is well known that the
family of $(1,1)$-knots contains all torus knots (trivially) and
all two-bridge knots (see \cite{KS}) in $\S^3$.

In this paper we develop an algebraic representation of
$(1,1)$-knots through elements of $MCG_2(T)$, the mapping class
group of the twice punctured torus. In Section \ref{(1,1) e MCG}
we establish the connection between the two objects. In Section
\ref{standard} we prove that every $(1,1)$-knot in a lens space
$L(p,q)$ can be represented by an element of $MCG_2(T)$ which is
the composition of an element of a certain rank two free subgroup
and of a standard element only depending on the ambient space
$L(p,q)$. This representation will be called ``standard''. As a
notable application, in Sections \ref{torus} and \ref{due ponti}
we obtain standard representations for the two most important
classes of $(1,1)$-knots in $\S^3$: the torus knots and the
two-bridge knots. Moreover, applying certain results obtained in
\cite{CM}, we give explicit cyclic presentations for the
fundamental groups of all cyclic branched coverings of torus knots
of type $(k,ck+2)$, with $c,k>0$ and $k$ odd.

In what follows, the symbol $L(p,q)$ will denote any lens space,
including $\S^3=L(1,0)$ and $\S^1\times\S^2=L(0,1)$. Moreover,
homotopy and homology classes will be denoted with the same symbol
of the representing loops.

\bigskip

\section{$\mathbf{(1,1)}$-knots and $\mathbf{MCG_2(T)}$}
\label{(1,1) e MCG} Let $F_g$  be a closed orientable surface of
genus $g$ and let $\mathcal{P}=\{P_1,\ldots ,P_n\}$ be a finite
set of distinguished points of $F_g$, called \textit{punctures}.
We denote by $\mathcal{H}(F_g,\mathcal{P})$ the group of
orientation-preserving homeomorphisms $h:F_g\to F_g$ such that
$h(\mathcal{P})=\mathcal{P}$. The \textit{punctured mapping class
group} of $F_g$ relative to $\mathcal{P}$ is the group of the
isotopy classes of elements of $\mathcal{H}(F_g,\mathcal{P})$. Up
to isomorphism, the punctured mapping class group of a fixed
surface $F_g$ relative to $\mathcal{P}$ only depends on the
cardinality $n$ of $\mathcal{P}$. Therefore, we can simply speak
of the {\it $n$-punctured mapping class group\/} of $F_g$,
denoting it by $MCG_n(F_g)$. Moreover, for isotopy classes we will
use the same symbol of the representing homeomorphisms.

The \textit{$n$-punctured pure mapping class group} of $F_g$ is
the subgroup $PMCG_n(F_g)$ of $MCG_n(F_g)$ consisting of the
elements pointwise fixing the punctures. There is a standard exact
sequence
$$1\to PMCG_n(F_g) \to MCG_n(F_g)\to \Sigma_n \to 1,$$
where $\Sigma_n$ is the symmetric group on $n$ elements. A
presentation of all punctured mapping class groups can be found in
\cite{Ge} and in \cite{LP}.

In this paper we are interested in the two-punctured mapping class
group of the torus $MCG_2(T)$. According to previously cited
papers, a set of generators for $MCG_2(T)$ is given by a rotation
$\r$ of $\p$ radians which exchanges the punctures and the
right-handed Dehn twists $\t_{\a},\t_{\b},\t_ {\bb}$ around the
curves $\a,\b,\bb$ respectively, as depicted in Figure \ref{Fig.
2}. Since $\r$ commutes with the other generators, we have
$$MCG_2(T)\cong PMCG_2(T)\oplus\Z_2.$$

\begin{figure}[ht]
\begin{center}
\includegraphics*[totalheight=4cm]{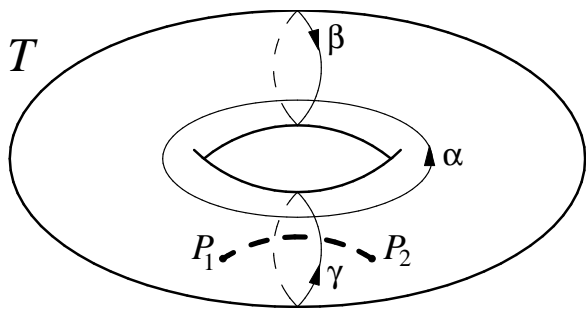}
\end{center}
\caption{Generators of $MCG_2(T)$.} \label{Fig. 2}
\end{figure}

The following presentation for $PMCG_2(T)$ has been obtained in
\cite{PS}:
\begin{equation}\label{presentation}
\langle\,\t_{\a},\t_{\b},\t_{\bb}\ | \
\t_{\a}\t_{\b}\t_{\a}=\t_{\b}\t_{\a}\t_{\b},\,
\t_{\a}\t_{\bb}\t_{\a}=\t_{\bb}\t_{\a}\t_{\bb},\,
\t_{\b}\t_{\bb}=\t_{\bb}\t_{\b},\,
(\t_{\a}\t_{\b}\t_{\bb})^{4}=1\,\rangle.
\end{equation}

The group $PMCG_2(T)$ (as well as $MCG_2(T)$) naturally maps by an
epimorphism to the mapping class group of the torus $MCG(T)\cong
SL(2,\Z)$, which is generated by $t_{\a}$ and $t_{\b}=t_{\bb}$. So
we have an epimorphism
$$\O:PMCG_2(T)\to SL(2,\Z)$$ defined by $\O(t_{\a})=
\left(\begin{matrix}1&0\\1&1\end{matrix}\right)$ and
$\O(t_{\b})=\O( t_{\bb})=
\left(\begin{matrix}1&-1\\0&1\end{matrix}\right)$.

\medskip

\begin{figure}[h]
\begin{center}
\includegraphics*[totalheight=7.5cm]{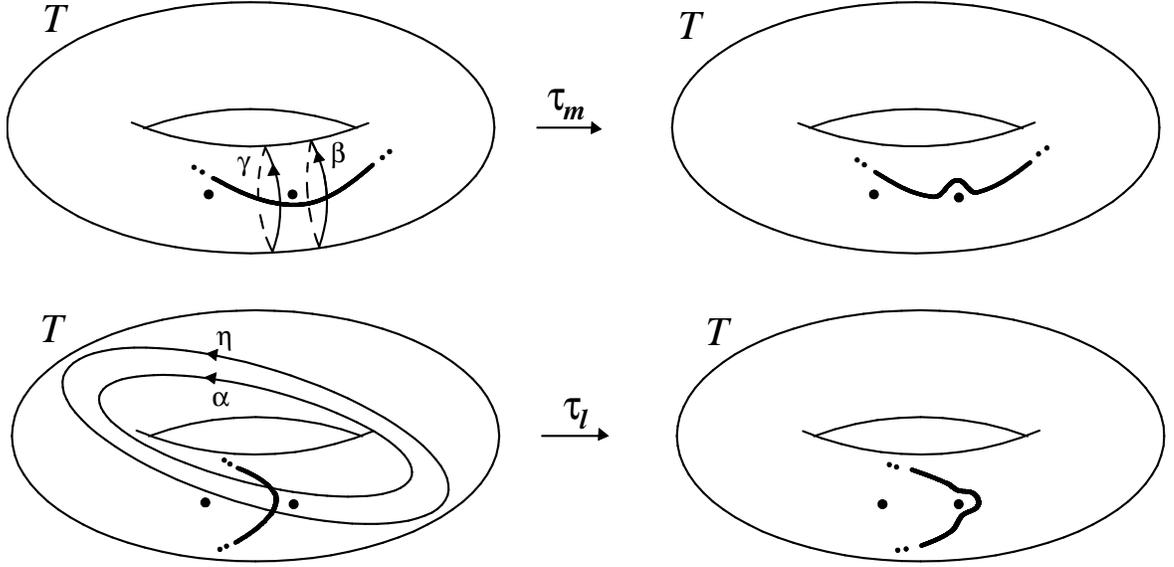}
\end{center}
\caption{Action of $\tau_m$ and $\tau_l$.} \label{Fig. 3}
\end{figure}

The group $\ker\O$ will play a fundamental role in our discussion.
In order to investigate its structure, let us consider the  two
elements $\tau_m=t_{\b}t_{\gamma}^{-1}$ and
$\tau_l=t_{\aa}t_{\a}^{-1}$, where $t_{\aa}$ is the right-handed
Dehn twist around the curve $\aa$ depicted in Figure \ref{Fig. 3}.
The effect of $\tau_m$ and $\tau_l$ is to slide one puncture (say
$P_2$) respectively along a meridian and along a longitude of the
torus, as shown in Figure~\ref{Fig. 3}. Observe that, since
$\aa=\tau_m^{-1}(\a)$,  we have $t_{\aa}=\tau_m^{-1}
t_{\a}\tau_m$.

The following result can be obtained from \cite[Th. 1]{Bi3} and
\cite[Th. 5]{Bi1} by classical techniques.

\begin{proposition} \label{nucleo} The group $\ker\O$ is freely generated by
$\tau_m=t_{\b}t_{\gamma}^{-1}$ and $\tau_l=t_{\aa}t_{\a}^{-1}$,
where $t_{\aa}=\tau_m^{-1} t_{\a}\tau_m$.
\end{proposition}

Now, let $K\subset L(p,q)$ be a $(1,1)$-knot with $(1,1)$-decomposition
$(L(p,q),K)=(H,A)\cup_{\ps}(H',A')$ and let $\mu :(H,A)\to
(H',A')$ be a fixed orientation-reversing homeomorphism, then
$\f=\ps\mu_{|\partial H}$ is an orientation-preserving
homeomorphism of $(\partial H,\partial A)=(T,\{P_1,P_2\})$.
Moreover, since two isotopic attaching homeomorphisms produce
equivalent $(1,1)$-knots, we have a natural surjective map from
the twice punctured mapping class group of the torus $MCG_2(T)$ to
the class $\K_{1,1}$ of all $(1,1)$-knots
$$\Theta:\f\in MCG_2(T)\mapsto K_{\f}\in \K_{1,1}.$$

If $\O(\f)=\left(\begin{matrix}q&s\\p&r\end{matrix}\right)$, then
$K_{\f}$ is a $(1,1)$-knot in the lens space $L(\vert p\vert,\vert
q\vert)$ \cite[p. 186]{BZ}, and therefore it is a knot in $\S^3$
if and only if $p=\pm 1$.

As will be proved in Section 3, we have the following ``trivial''
examples:
\begin{itemize}
\item[i)] if either $\f=1$ or $\f=t_{\b}$ or $\f=t_{\bb}$, then
$K_{\f}$ is the trivial knot in $\S^1\times\S^2$;
\item[ii)] if $\f=t_{\a}$, then $K_{\f}$ is the trivial knot in $\S^3$.
\end{itemize}

Moreover, it is possible to prove that if
$\f=t_{\a}t_{\b}t_{\a}t_{\a}t_{\gamma}t_{\a}$, then $K_{\f}$ is
the knot $\S^1\times\{P\}\subset\S^1\times\S^2$, where $P$ is any
point of $\S^2$. So, in this case, $K_{\f}$ is a standard
generator for the first homology group of $\S^1\times\S^2$.

Every element $\f$ of $MCG_2(T)$ can be written as $\f=\f'\r^k$,
$k\in\{0,1\}$, where $\f'\in PMCG_2(T)$. Since $\r$ can be
extended to a homeomorphism of the pair $(H,A)$, the $(1,1)$-knots
$K_{\f}$ and $K_{\f'}$ are equivalent. So,  for our discussion it
is enough to consider the restriction
$$\Theta'=\Theta_{|PMCG_2(T)}:\f\in PMCG_2(T)\mapsto K_{\f}\in \K_{1,1}.$$

\bigskip

\section{Standard decomposition}
\label{standard} In this section we show that every $(1,1)$-knot
$K\subset L(p,q)$ admits a representation by the composition of an
element in $\ker\O$ and an element which only depends on $L(p,q)$.
A representation of this type will be called ``standard''. Note
that a similar result, using a rank three free subgroup of
$MCG_2(T)$, has been obtained in \cite[Theorem 3]{CK}.

First of all, we deal with trivial knots in lens space. Let $\T$
be the subgroup of $PMCG_2(T)$ generated by $t_{\a}$ and $t_{\b}$.
There exists a disk $D\subset H$, with $A\cap D=A\cap\partial D=A$
and $\partial D-A\subset T$, such that
$D\cap\a=D\cap\b=\emptyset$. So any element of $\T$ produces a
trivial knot in a certain lens space. On the other hand, any
trivial knot in a lens space admits a representation through an
element of $\T$, as will be proved in Proposition
\ref{trivialknots}.

We need a preparatory result.

\begin{lemma} \label{matriciequivalenti}
Let $K$ be a $(1,1)$-knot in $L(p,q)$. Then, for each $r,s\in\Z$
such that $qr-ps=1$ there exists $\f\in PMCG_2(T)$, with
$\O(\f)=\left(\begin{matrix}q&s\\p&r\end{matrix}\right)$, such
that $K=K_{\f}$.
\end{lemma}
\begin{proof} Let $K=K_{\bar\f}$, with $\O(\bar\f)=\left(\begin{matrix}q&\bar s\\p&\bar
r\end{matrix}\right)$. Since $q\bar r-p\bar s=1$, there exist
$c\in\Z$ such that $r=\bar r+cp$ and $s=\bar s+cq$. If $\f=\bar\f
t_{\b}^{-c}$, we have $K_{\f}=K_{\bar\f}$, since $t_{\b}^{-c}$ can
be extended to a homeomorphism of the pair $(H,A)$. Moreover
$\O(\f)=\O(\bar\f)\O(t_{\b}^{-c})=\left(\begin{matrix}q&\bar
s\\p&\bar r\end{matrix}\right)\ \left(\begin{matrix}1&c\\0&1
\end{matrix}\right)
=\left(\begin{matrix}q&\bar s+cq\\p&\bar r+cp\end{matrix}\right).$
\end{proof}

\medskip

For integers $p,q$ such that $0<q<p$ and $\gcd(p,q)=1$ consider
the sequence of equations of the Euclidean algorithm (with
$r_0=p$, $r_1=q$):
\begin{eqnarray*} \label{nodobanale} r_0&=&a_1 r_1+r_2\\ r_1&=&a_2 r_2+r_3\\
\vdots&\ &\\ r_{m-2}&=&a_{m-1} r_{m-1}+r_m \\ r_{m-1}&=&a_m r_m,
\end{eqnarray*} with $r_1>r_2>\cdots>r_{m-1}>r_m=1$.

The $a_i$'s are the coefficients of the continued fraction
$$\frac {p}{q}  = a_1+ \frac{1}{a_2 +\frac{1} {a_3 + \cdots +
\frac{1}{a_m}}}\,.$$ In the following we will use the notation
$p/q=[a_1,a_2,\ldots,a_m]$.

\begin{proposition} \label{trivialknots}
\begin{itemize}
\item The trivial knot in $\S^3=L(1,0)$ is represented by
$\f_{1,0}=\t_{\b}\t_{\a}t_{\b}$. \item The trivial knot in
$\S^1\times\S^2=L(0,1)$ is represented by $\f_{0,1}=1$. \item Let
$p,q$ be integers such that $0<q<p$ and $\gcd(p,q)=1$. If
$p/q=[a_1,a_2,\ldots,a_m]$, then the trivial knot in the lens
space $L(p,q)$ is represented by
$$\f_{p,q}=\begin{cases}
t_{\a}^{a_1}t_{\b}^{-a_2}\cdots t_{\a}^{a_m}& \text{ if $m$ is odd}\\
t_{\a}^{a_1}t_{\b}^{-a_2}\cdots t_{\b}^{-a_m}t_{\b}t_{\a}t_{\b} &
\text{if $m$ is even}\end{cases}.$$
\end{itemize}
\end{proposition}
\begin{proof}
Since all the involved homeomorphisms belong to $\T$, all the
knots are trivial. It is easy to check (see also \cite[p.
186]{BZ}) that, for suitable $r,s\in\Z$, we have:
$$
\left(\begin{matrix}q&s\\p&r\end{matrix}\right)=\begin{cases}\left(\begin{matrix}1&0\\a_1&1\end{matrix}\right)\
\left(\begin{matrix}1&a_2\\0&1\end{matrix}\right)\cdots
\left(\begin{matrix}1&0\\a_m&1\end{matrix}\right)&\text{if  $m$ is
odd},\\ \\ \left(\begin{matrix}1&0\\a_1&1\end{matrix}\right)\
\left(\begin{matrix}1&a_2\\0&1\end{matrix}\right)\cdots
\left(\begin{matrix}1&a_m\\0&1\end{matrix}\right)\
\left(\begin{matrix}0 & - 1\\
1&0\end{matrix}\right)&\text{ if  $m$  is even}.
\end{cases}$$
Since
$\O(\t_{\a}^{a_i})=\left(\begin{matrix}1&0\\a_i&1\end{matrix}\right)$,
$\O(\t_{\b}^{a_i})=\left(\begin{matrix}1&-a_i\\0&1\end{matrix}\right)$,
and $\O(\t_{\b}\t_{\a}t_{\b})=\left(\begin{matrix}0 & -1\\
1&0\end{matrix}\right)$, the statement is obtained.
\end{proof}

\medskip

Now we can prove the result announced at the beginning of this
section.

\begin{theorem} \label{decomposition} Let $K$ be a $(1,1)$-knot
in $L(p,q)$. Then there exist $\f',\f''\in\ker\O$ such that
$K=K_{\f}$, with $\f=\f'\f_{p,q}=\f_{p,q}\f''$.
\end{theorem}
\begin{proof} By Lemma \ref{matriciequivalenti}, there
exists $\f$, with $\O(\f)=\O(\f_{p,q})$, such that $K=K_{\f}$. It
suffices to define $\f'=\f\f_{p,q}^{-1}$ and
$\f''=\f_{p,q}^{-1}\f$.
\end{proof}

\medskip

A representation $\f\in PMCG_2(T)$ of a $(1,1)$-knot will be
called {\it standard\/} if $\f$ is of the type
described in the previous theorem.

We point out that $(1,1)$-knots admit different (usually
infinitely many) standard representations. For example $\tau_m^c$
represents the trivial knot in \hbox{$\S^1\times\S^2$}, for all
$c\in\Z$.

\section{Representation of torus knots}
\label{torus}
In this section we give a standard representation for all torus
knots in $\S^3$. Let $K={\bf t}(k,h)$ be a torus knot of type
$(k,h)$. Then $\gcd(k,h)=1$, and we can assume that $K$ lies on
the boundary $T=\partial H$ of a genus one handlebody $H$
canonically embedded in $\S^3$. The homology class of $K$ is $hl +
km$, where $l$ and $m$ respectively denote a longitude and a
meridian of $T$. By slightly pushing (the interior of) an arc
$A'\subset K$ outside $H$ and $K-A'$ inside $H$, we obtain an
obvious $(1,1)$-decomposition of $K$. Observe that $0<\vert
k\vert<h$ can be assumed without loss of generality (see \cite[p.
45]{BZ}).

\begin{figure}[ht]
\begin{center}
\includegraphics*[totalheight=4.5cm]{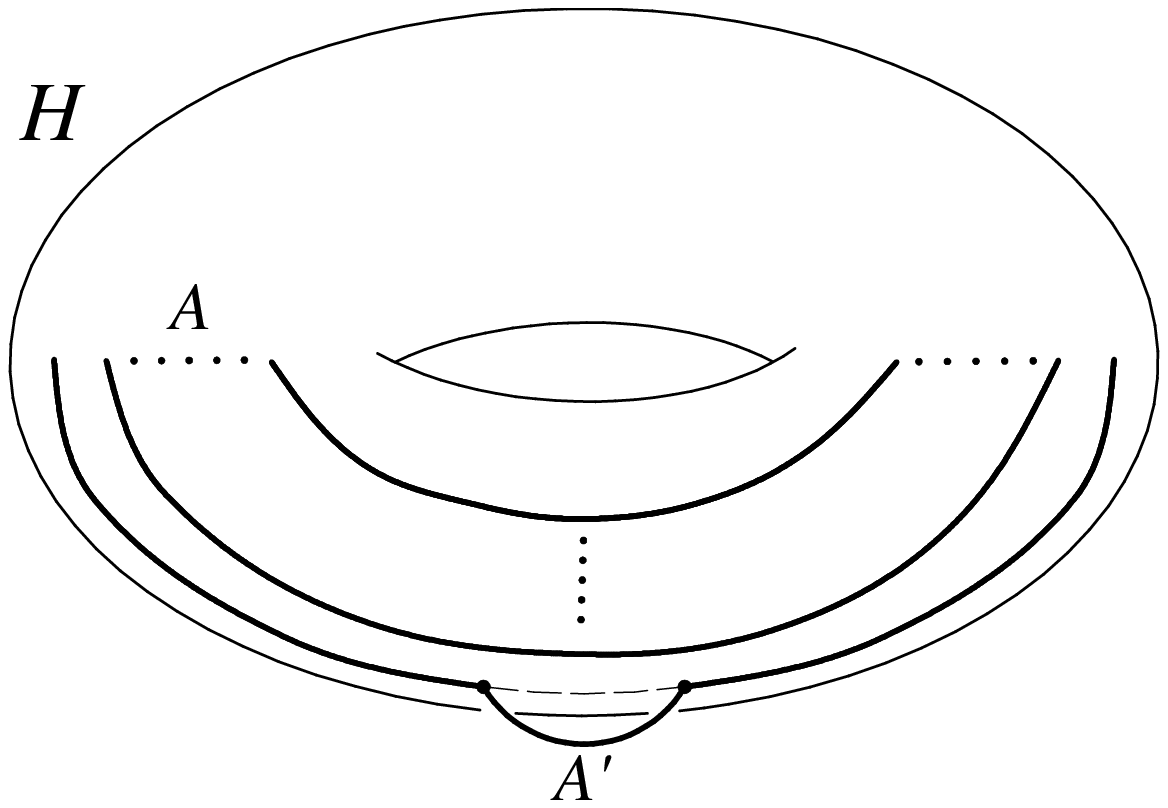}
\end{center}
\caption{} \label{Fig. 5}
\end{figure}

In the next statement $\lfloor x \rfloor$ denotes the integral
part of $x$.

\begin{theorem} \label{torusknots} The torus knot ${\bf
t}(k,h)\subset\S^3$ is the $(1,1)$-knot $K_{\f}$ with:
$$\f=\prod_{i=1}^{h}(\tau_m^{\lfloor (i-1)k/h
\rfloor-\lfloor ik/h \rfloor}\tau_l^{-1})\t_{\b}\t_{\a}t_{\b},$$
where $\tau_m=t_{\b}t_{\gamma}^{-1}$ and
$\tau_l=\tau_m^{-1}t_{\a}\tau_m t_{\a}^{-1}$.
\end{theorem}

\begin{figure}[ht]
\begin{center}
\includegraphics*[totalheight=8cm]{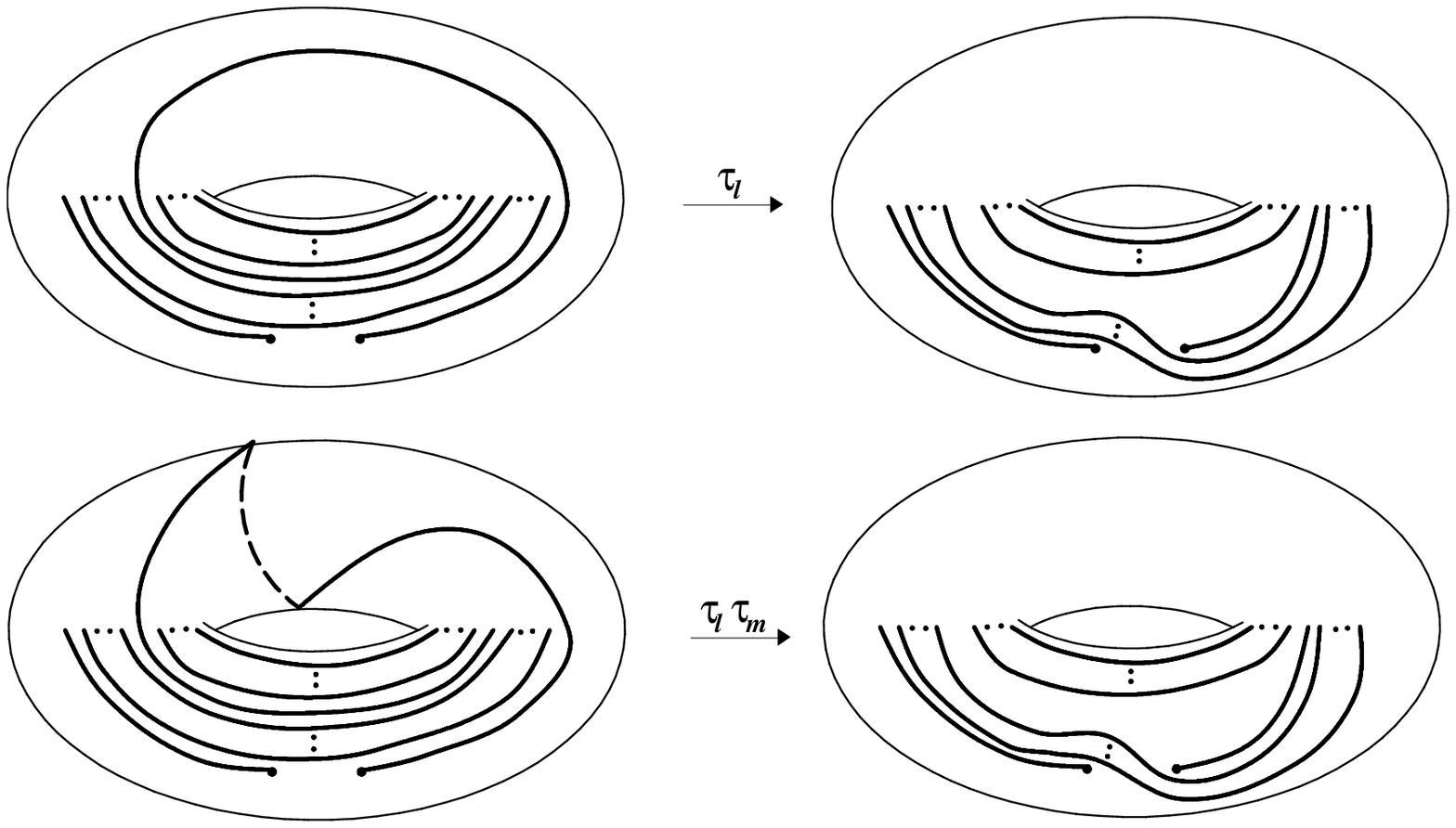}
\end{center}
\caption{} \label{Fig. 6}
\end{figure}

\begin{proof} Up to isotopy, we can suppose that the arc
$A=K_{\f}-\textup{int}(A')$ lies on $\partial H$, as in Figure
\ref{Fig. 5}. The arc $A$ can be transformed into an arc
$\tilde{A}$ in such a way that $\tilde{A}\cup A'$ is a trivial
knot in $\S^3$, represented by the standard homeomorphism
$\f_{1,0}=\t_{\b}\t_{\a}t_{\b}$, via a suitable sequence of
homeomorphisms $\tau_l$ and $\tau_m$, according to the following
algorithm. Consider the sequence of equations:
\begin{eqnarray*}
& &k=q_1h+r_1,\\
& &2k=q_2h+r_2,\\
& &\vdots\\
& &hk=q_{h}h+r_{h},
\end{eqnarray*}

\begin{figure}
\begin{center}
\includegraphics*[totalheight=18cm]{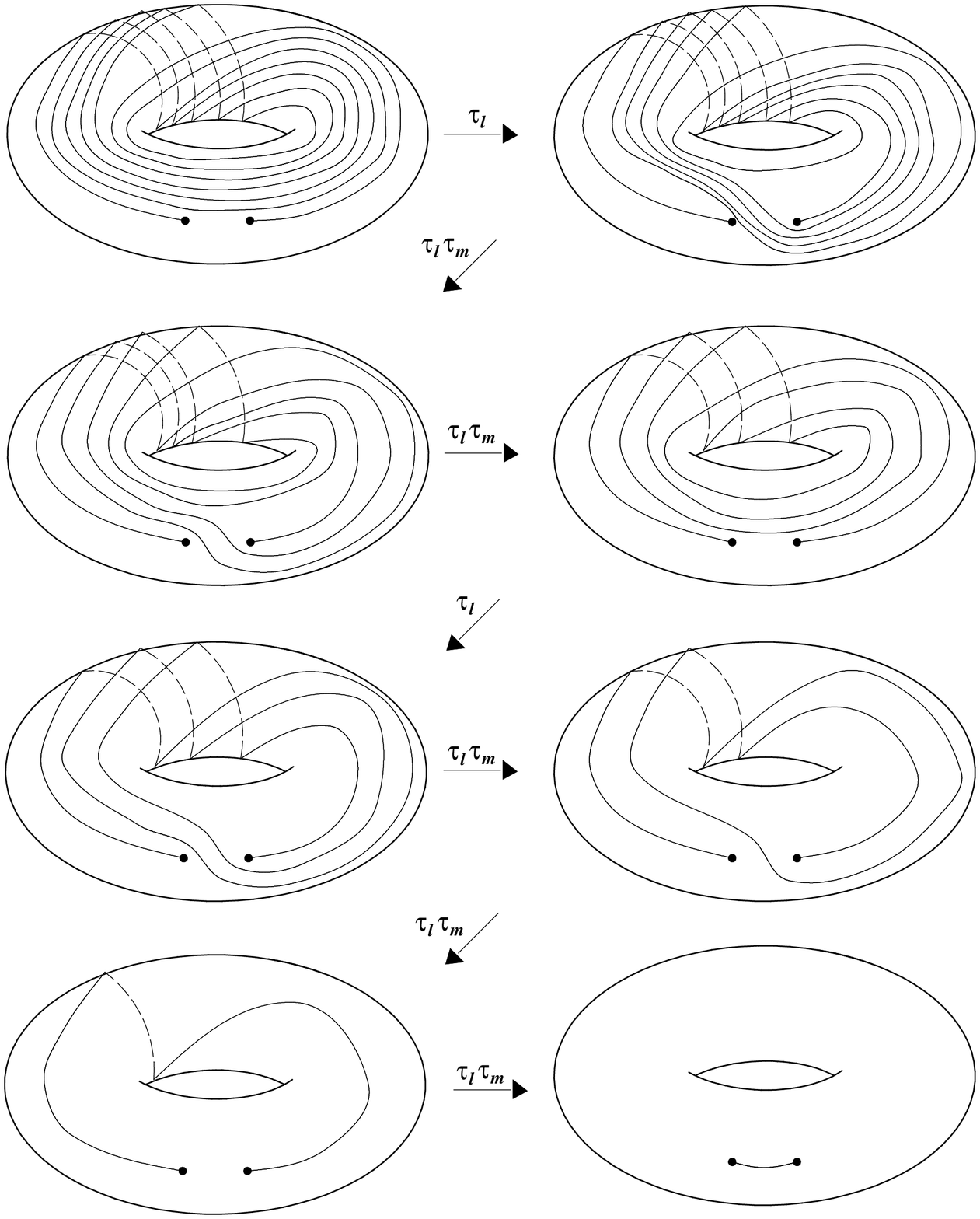}
\end{center}
\caption{Trivialization of $\mathbf{t}(5,7)$.} \label{Fig. 7}
\end{figure}

\noindent where $0\le r_i<h$, for $i=1,\ldots h$. Moreover, define $q_0=0$.
So $q_i=\lfloor ik/h \rfloor$, for $i=0,1,\ldots h$. Now define
the homeomorphisms $\f_i=\tau_l\tau_m^{q_i-q_{i-1}}$, for
$i=1,\ldots,h$. Figure \ref{Fig. 6} depicts the effect of $\tau_l$
and $\tau_l\tau_m$ on $A$. As a consequence, the homeomorphism
$\phi=\f_h\f_{h-1}\cdots\f_1$ transforms the arc $A$ into the arc
$\tilde{A}$ (Figure 7 shows the case ${\bf t}(5,7)$), and
therefore we have $\f_{1,0}=\phi\f$. So $\phi^{-1}\f_{1,0}$
represents the torus knot ${\bf t}(k,h)$.
\end{proof}

For example, ${\bf t}(5,7)=K_{\f}$, with
$\f=\tau_l^{-1}(\tau_m^{-1}\tau_l^{-1})^2\tau_l^{-1}(\tau_m^{-1}\tau_l^{-1})^3\t_{\b}\t_{\a}t_{\b}$
(see Figure \ref{Fig. 7}).

\medskip

As a consequence, we obtain a cyclic presentation for the
fundamental group for all cyclic branched coverings of a
particular class of torus knots.

\begin{proposition} \label{torus2}
The fundamental group of the $n$-fold cyclic branched covering of
the torus knot $\mathbf{t}(k,ck+2)$, with $k>1$ odd and $c>0$,
admits the cyclic presentation $G_n(w)$, where $w$ is equal to
$$\prod_{i=0}^{(k-3)/2}
(\prod_{j=0}^{c(k-1)/2}x_{1-i(ck+2)+jk}
\prod_{l=0}^{c(k+1)/2}x_{ck(k-1)/2-i(ck+2)-lk}^{-1})
 \prod_{m=0}^{c(k-1)/2}x_{1-(k-1)(ck+2)/2+mk}$$
(subscripts are taken modulo $n$).
\end{proposition}
\begin{proof} Let $r=(k-1)/2$. From Theorem \ref{torusknots} we have $\mathbf{t}(k,ck+2)=K_{\f}$ with
$\f=(\tau_l^{-c}\tau_m^{-1})^r\tau_l^{-1}(\tau_l^{-c}\tau_m^{-1})^r\tau_l^{-c}\tau_m^{-1}\tau_l^{-1}t_{\beta}t_{\alpha}t_{\beta}$.
Applying \cite[Proposition 1]{CM}, we obtain
$\p_1(\S^3-\mathbf{t}(k,ck+2))=\langle\la,\lbb\,|\,r(\la,\lbb)\rangle,$
with
$r(\la,\lbb)=(\lbb^{-1}\la^{cr+1}\lbb^{-1}\la^{-c(r+1)-1})^r\lbb^{-1}\la^{cr+1}$.
Then
$H_1(\S^3-\mathbf{t}(k,ck+2))=\langle\a,\bb\,|\,\a-k\bb\rangle$.
Since, up to equivalence, $\o_f(\bb)=1$,  we have $\o_f(\a)=k$. We
set $\la= x\lbb^k$, therefore
$\p_1(\S^3-\mathbf{t}(k,ck+2))=\langle x ,\lbb\,|\,\bar r(
x,\lbb)\rangle,$ with $\bar r( x ,\lbb)=(\lbb^{-1}(
x\lbb^{k})^{1+c(k-1)/2}\lbb^{-1}( \lbb^{-k} x^{-1}
)^{1+c(k+1)/2})^{(k-1)/2}\lbb^{-1}( x\lbb^{k})^{1+c(k-1)/2}$. The
statement derives from a straightforward application of
\cite[Theorem 7]{CM},
\end{proof}

\medskip

For example, the fundamental group of the $n$-fold cyclic branched
covering of ${\bf t}(5,7)$ admits the cyclic presentation
$G_n(w)$, where $$w=
x_{15}x_{20}x_{25}x_{24}^{-1}x_{19}^{-1}x_{14}^{-1}x_{9}^{-1}x_{8}x_{13}x_{18}
x_{17}^{-1}x_{12}^{-1}x_{7}^{-1}x_{2}^{-1}x_{1}x_{6}x_{11}.$$

\bigskip

\section{Representation of two-bridge knots}
\label{due ponti}
\begin{figure}[ht]
\begin{center}
\includegraphics*[totalheight=1.8cm]{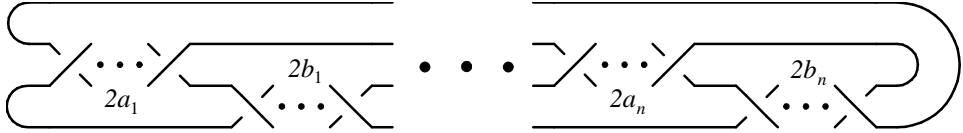}
\end{center}
\caption{Conway presentation for two-bridge knots.}
\label{figconway}
\end{figure}

In this section we give a standard representation for all
two-bridge knots in $\S^3$. Let ${\bf b}(a/b)$ be a non-trivial
two-bridge knot in $\S^3$ of type $(a,b)$. Then we can assume
$\gcd(a,b)=1$, $a$ odd, $b$ even and $0<\vert b\vert<a$, without
loss of generality (see \cite[Ch. 12B]{BZ}). It is known that
${\bf b}(a/b)$ admits a Conway presentation with an even number of
even parameters $[2a_1,2b_1,\ldots ,2a_n,2b_n]$ (see Figure
\ref{figconway}), satisfying the following relation:
$$\frac {a}{b}  = 2a_1+ \frac{1}{2b_1 +\frac{1} {2a_2 + \cdots +\frac{1}{2b_n}}}\,.$$

\begin{figure}
\begin{center}
\includegraphics*[totalheight=18cm]{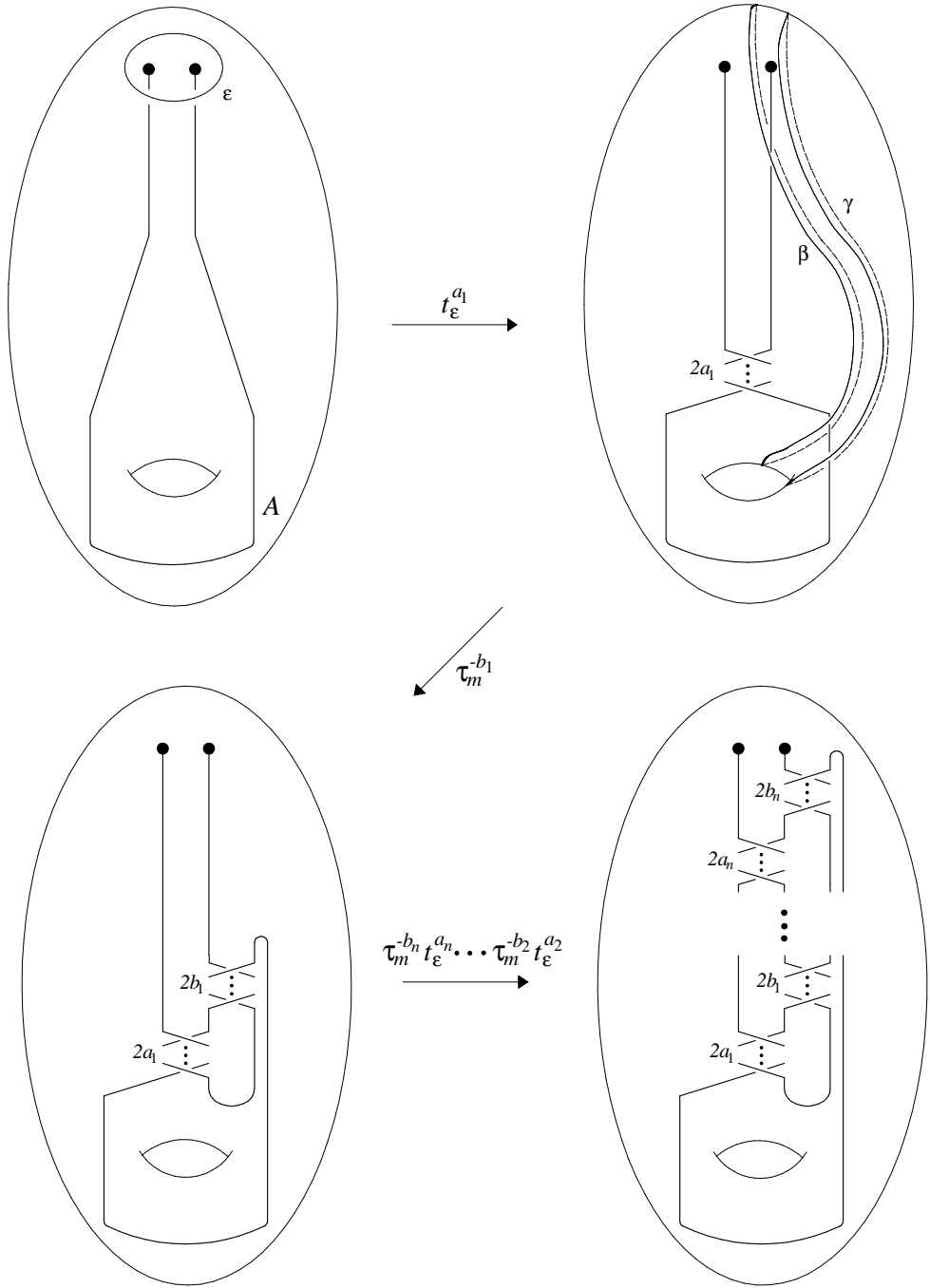}
\end{center}
\caption{Standard representation of two-bridge knots.}
\label{figtwobridge}
\end{figure}

\begin{theorem} \label{twobridge}
The two-bridge knot ${\bf b}(a/b)\subset\S^3$ having Conway
parameters $[2a_1,2b_1,\ldots ,2a_n,2b_n]$ is the $(1,1)$-knot
$K_{\f}$ with:
$$\f=\t_{\b}\t_{\a}t_{\b}\tau_m^{-b_n}t_{\e}^{a_n}\cdots
\tau_m^{-b_1}t_{\e}^{a_1},$$ where
$t_{\e}=\tau_l^{-1}\tau_m\tau_l\tau_m^{-1}$ is the right-handed
Dehn twist around the curve $\e$ depicted in Figure
\ref{figtwobridge}.
\end{theorem}

\begin{figure}
\begin{center}
\includegraphics*[totalheight=4 cm]{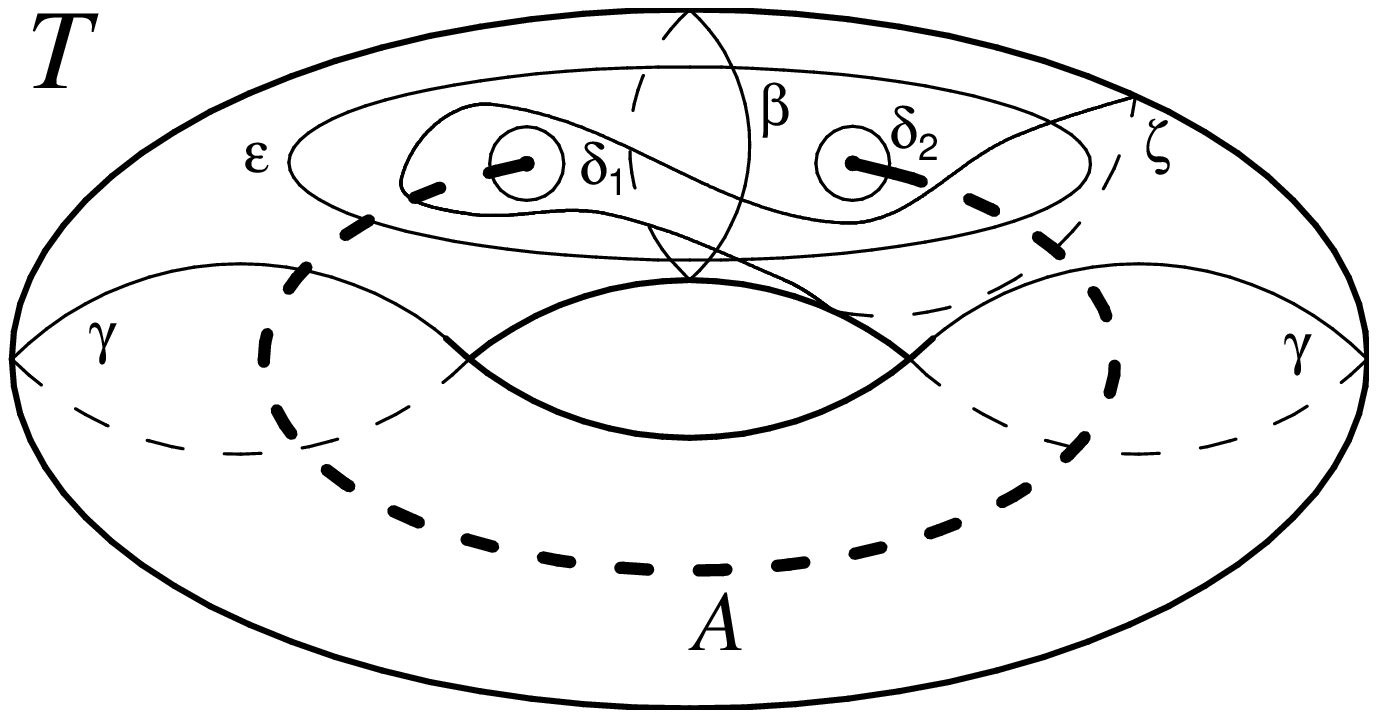}
\end{center}
\caption{} \label{Fig. 17}
\end{figure}

\begin{proof}
Figure \ref{figtwobridge} shows the result of the application of
$\tau_m^{-b_n}t_{\e}^{a_n}\cdots \tau_m^{-b_1}t_{\e}^{a_1}$. By
applying $\f_{1,0}=\t_{\b}\t_{\a}t_{\b}$ we obtain the two-bridge
knot with Conway parameters $[2a_1,2b_1,\ldots ,2a_n,2b_n]$.

Now we show that $t_{\e}=\tau_l^{-1}\tau_m\tau_l\tau_m^{-1}$ (note
that no disk bounded by $\e$ and properly embedded in $H$ is
disjoint from $A$). Referring to  Figure \ref{Fig. 17}, the
following ``lantern'' relation
$\t_{\bb}^2\t_{\delta_1}\t_{\delta_2}=\t_{\e}\t_{\b}\t_{\zeta}$
holds (see \cite{W}). So we obtain
$\zeta=\t_{\a}\t_{\bb}\t_{\b}^{-1}\t_{\a}^{-1}(\bb)$ and therefore
$\t_{\zeta}=\t_{\a}\t_{\bb}\t_{\b}^{-1}\t_{\a}^{-1}\t_{\bb}\t_{\a}\t_{\b}\t_{\bb}^{-1}\t_{\a}^{-1}$.
Since $\t_{\delta_1}=\t_{\delta_2}=1$ we have
$\t_{\e}=\t_{\bb}^2\t_{\zeta}^{-1}\t_{\b}^{-1}=\t_{\bb}^2\t_{\a}\t_{\bb}\t_{\b}^{-1}\t_{\a}^{-1}\t_{\bb}^{-1}\t_{\a}\t_{\b}\t_{\bb}^{-1}
\t_{\a}^{-1}\t_{\b}^{-1}.$ Now, using the relations of
(\ref{presentation}) we get
\begin{eqnarray*}\t_{\e} &=&\t_{\bb}^2\t_{\a}\t_{\bb}\t_{\b}^{-1}\t_{\a}^{-1}\t_{\bb}^{-1}\t_{\a}\t_{\b}\t_{\bb}^{-1}
\t_{\a}^{-1}\t_{\b}^{-1}=\t_{\bb}\t_{\a}\t_{\bb}\t_{\a}\t_{\b}^{-1}\t_{\a}^{-1}\t_{\bb}^{-1}\t_{\a}\t_{\b}\t_{\bb}^{-1}
\t_{\a}^{-1}\t_{\b}^{-1}=\\
\ &=&\t_{\bb}\t_{\a}\t_{\bb}\t_{\b}^{-1}\t_{\a}^{-1}\t_{\b}\t_{\bb}^{-1}\t_{\a}\t_{\b}\t_{\bb}^{-1}
\t_{\a}^{-1}\t_{\b}^{-1}=\t_{\bb}\t_{\a}\t_{\b}^{-1}\t_{\bb}\t_{\a}^{-1}\t_{\bb}^{-1}\t_{\b}\t_{\a}\t_{\b}\t_{\bb}^{-1}
\t_{\a}^{-1}\t_{\b}^{-1}=\\
\ &=&\t_{\bb}\t_{\a}\t_{\b}^{-1}\t_{\a}^{-1}\t_{\bb}^{-1}\t_{\a}\t_{\b}\t_{\a}\t_{\b}\t_{\bb}^{-1}
\t_{\a}^{-1}\t_{\b}^{-1}=\t_{\bb}\t_{\b}^{-1}\t_{\a}^{-1}\t_{\b}\t_{\bb}^{-1}\t_{\a}\t_{\b}\t_{\a}\t_{\b}\t_{\bb}^{-1}
\t_{\a}^{-1}\t_{\b}^{-1}=\\
\ &=&\t_{\bb}\t_{\b}^{-1}\t_{\a}^{-1}\t_{\b}\t_{\bb}^{-1}\t_{\a}\t_{\a}\t_{\b}\t_{\a}\t_{\bb}^{-1}
\t_{\a}^{-1}\t_{\b}^{-1}=\t_{\bb}\t_{\b}^{-1}\t_{\a}^{-1}\t_{\b}\t_{\bb}^{-1}
\t_{\a}\t_{\a}\t_{\b}\t_{\bb}^{-1}\t_{\a}^{-1}\t_{\bb}\t_{\b}^{-1}=
\\
\ &=&\tau_m^{-1}\t_{\a}^{-1}\tau_m\t_{\a}\t_{\a}\tau_m\t_{\a}^{-1}\tau_m^{-1}=
\tau_m^{-1}\t_{\a}^{-1}\tau_m\t_{\a}\tau_m\tau_m^{-1}\t_{\a}\tau_m\t_{\a}^{-1}\tau_m^{-1}=\\
\
&=&\t_{\aa}^{-1}\t_{\a}\tau_m\t_{\aa}\t_{\a}^{-1}\tau_m^{-1}=\tau_l^{-1}\tau_m\tau_l\tau_m^{-1}.
\end{eqnarray*}
\end{proof}

\medskip

For example, the figure-eight knot ${\bf b}(5/2)$, which has
Conway parameters $[2,2]$, is the knot $K_{\f}$ with $\f
=\t_{\b}\t_{\a}\t_{\b}\tau_m^{-1}\t_{\e}$ (see Figure \ref{Fig.
19}).

\begin{figure}
\begin{center}
\includegraphics*[totalheight=17 cm]{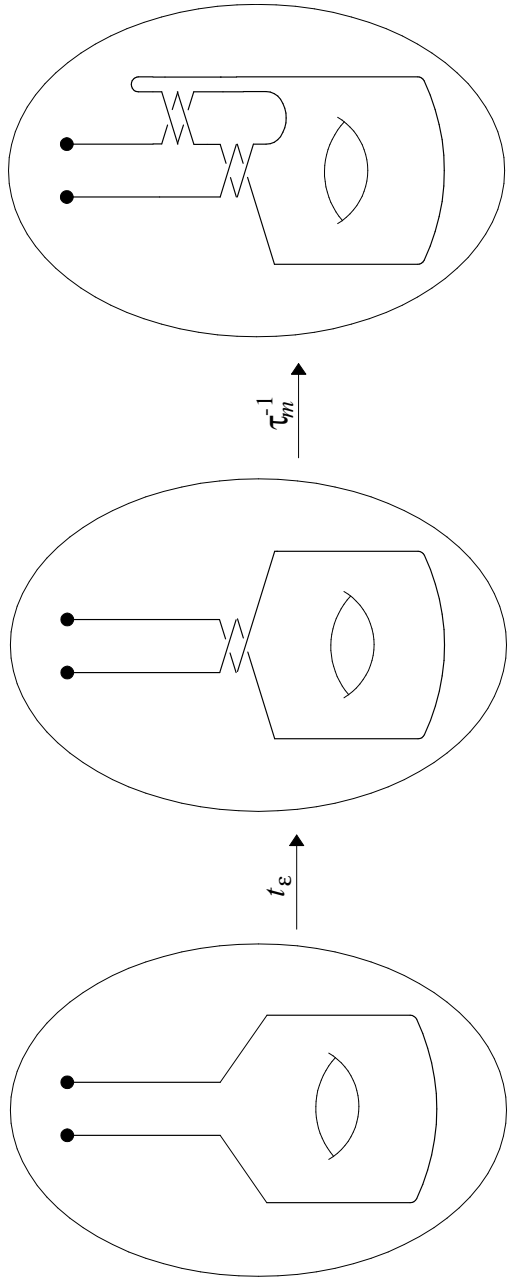}
\end{center}
\caption{Standard representation of the figure-eight knot.}
\label{Fig. 19}
\end{figure}

\newpage


\bigskip\bigskip

\noindent \textit{Acknowledgements.} The authors  would like to
thank Sylvain Gervais and Andrei Vesnin for their helpful
suggestions. We also would like to thank the Referee for his
valuable comments and remarks. Work performed under the auspices
of the G.N.S.A.G.A. of I.N.d.A.M. (Italy) and the University of
Bologna, funds for selected research topics.

\bigskip


\bigskip

\vspace{15 pt} {ALESSIA CATTABRIGA, Department of Mathematics,
University of Bologna, Italy. E-mail: cattabri@dm.unibo.it}

\vspace{15 pt} {MICHELE MULAZZANI, Department of Mathematics and
C.I.R.A.M., University of Bologna, Italy. E-mail:
mulazza@dm.unibo.it}

\end{document}